\documentclass[12pt,thmsa]{article}
\usepackage{amssymb}
\usepackage{amsfonts}
\usepackage{sw20bams}

\input tcilatex
\begin{document}

\author{Steven R. Finch}
\title{Mutually Equidistant Spheres that Intersect}
\date{January 23, 2013}
\maketitle

\begin{abstract}
The setting for this brief paper is $\mathbb{R}^{3}$. \ Distance between two
spheres is understood as distance $\delta $ between spherical centers. \ For
instance,\ a Reuleaux tetrahedron $T$ is the intersection of four unit balls
satisfying $\delta =1$ pairwise. Volume and surface area of $T$ are already
well-known; our humble contribution is to calculate the mean width of $T$.
\end{abstract}

\footnotetext{%
Copyright \copyright\ 2013 by Steven R. Finch. All rights reserved.}Two
earlier papers \cite{Fi1, Fi2} were devoted to convex hulls involving disks.
We will here discuss intersections of balls. \ It is natural to compute
volume and surface area of such intersections. \ A\ third quantity, mean
width, is \textquotedblleft a new measure on three-dimensional solids that
enjoys equal rights with volume and surface area\textquotedblright\ \cite{Rt}%
. While results for $VL$ and $AR$ have appeared many times in the past, some
of our expressions for $MW$ may be new.

\section{Dihedron (Symmetric Lens)}

Consider, for simplicity, two unit spheres in $\mathbb{R}^{3}$ that pass
through each other's centers.\ \ The region enclosed by both spheres is
called a \textbf{spherical dihedron}, using language of \cite{G0, G3, G1, G2}%
. It is more commonly called a \textbf{symmetric lens}. Starting with spheres%
\[
\begin{array}{ccc}
\left( x-\dfrac{1}{2}\right) ^{2}+y^{2}+z^{2}=1, &  & \left( x+\dfrac{1}{2}%
\right) ^{2}+y^{2}+z^{2}=1%
\end{array}%
\]%
and defining%
\[
\begin{array}{ccc}
f(x,y)=\sqrt{1-\left( x+\dfrac{1}{2}\right) ^{2}-y^{2}}, &  & a(x)=\sqrt{%
1-\left( x+\dfrac{1}{2}\right) ^{2}}%
\end{array}%
\]%
we obtain%
\[
VL=8\dint\limits_{0}^{1/2}\;\;\dint\limits_{0}^{a(x)}f(x,y)dy\,dx=\frac{5\pi 
}{12}, 
\]%
\[
AR=8\dint\limits_{0}^{1/2}\;\;\dint\limits_{0}^{a(x)}\sqrt{%
1+f_{x}^{2}+f_{y}^{2}}\,dy\,dx=2\pi . 
\]%
Interiors of both faces have mean curvature $1$. The (unique) edge is a
circle with radius $f(0,0)=\sqrt{3}/2$, hence it has circumference $\sqrt{3}%
\pi $. \ The dihedral angle at the edge is $\pi /3$. By the
\textquotedblleft indirect approach\textquotedblright\ in \cite{Fi1}, we have%
\[
MW=\frac{1}{2\pi }AR+\frac{1}{4\pi }\frac{\pi }{3}\sqrt{3}\pi =1+\frac{\pi }{%
4\sqrt{3}}\text{.} 
\]

An obvious generalization allows $0<\delta <2$, where distance $\delta $
between spherical centers was taken to be $1$ above. \ Substantial work \cite%
{S1, S2, S3, S4, S5, S6, S7, S8, S9, S0} has been performed on intersections
of arbitrary collections of spheres in $\mathbb{R}^{3}$. We merely quote an
extended volume result \cite{He} for two unit spheres:%
\[
VL=\frac{4\pi }{3}\left( 1-\frac{3}{4}\delta +\frac{1}{16}\delta ^{3}\right)
. 
\]%
A\ symmetric lens is the union of two spherical caps. \ A different
parametrization makes use of the \textbf{angular radius} $\varphi $ of
either cap (also called the \textbf{angle of aperture} or \textbf{colatitude
angle} of the cap). \ From $\delta =2\cos (\varphi )$, it follows that \cite%
{Hw1} 
\[
VL=\frac{2\pi }{3}\left( 2-3\cos (\varphi )+\cos (\varphi )^{3}\right) , 
\]%
\[
AR=4\pi \left( 1-\cos (\varphi )\right) , 
\]%
\[
MW=2-2\cos (\varphi )+\left( \frac{\pi }{2}-\varphi \right) \sin (\varphi ). 
\]%
An unfortunate ambiguity appears in Figure 5f of \cite{Hw1}: the line
segment is intended to be normal to the rightmost spherical cap (\textit{not}
tangent to the leftmost spherical cap). Thus the indicated angle $\alpha $
is indeed equal to $\varphi $ (and \textit{not} equal to the \textbf{contact
angle} $\pi /2-\varphi $). Confusion occurs because $\alpha $ is chosen
close to $\pi /4$ in Figure 5f; the center picture in Figure 3 of \cite{Hw2}
helps to clarify matters.

We mention another error in the literature. The mean width in \cite{WL},
although defined to be $2\cdot MW$ (often encountered), gives results
inconsistent with \cite{Hw1}. \ See \cite{Oe} for another treatment of
overlapping spheres.

Higher $n$-dimensional results for both $VL$ and $AR$ can be inferred from 
\cite{T1, T2, T3, T4} that involve the Gauss hypergeometric function:%
\[
VL=\frac{\pi ^{n/2}}{\Gamma (1+n/2)}\left[ 1-\frac{2\,\Gamma (1+n/2)}{\sqrt{%
\pi }\,\Gamma ((n+1)/2)}\,_{2}F_{1}\left( \frac{1}{2},\frac{1-n}{2},\frac{3}{%
2},\cos (\varphi )^{2}\right) \cos (\varphi )\right] ,
\]%
\[
AR=\frac{2\pi ^{n/2}}{\Gamma (n/2)}\left[ 1-\frac{2\,\Gamma (n/2)}{\sqrt{\pi 
}\,\Gamma ((n-1)/2)}\,_{2}F_{1}\left( \frac{1}{2},\frac{3-n}{2},\frac{3}{2}%
,\cos (\varphi )^{2}\right) \cos (\varphi )\right] .
\]%
We have not attempted to find an analogous formula for $MW$. \ Older
alternative expressions for $VL$ and $AR$ are found by unraveling the
quermassintegrals $W_{0}$ and $k\,W_{1}$ in \cite{Hw2} (formula 69). \ New
and old results agree. \ Also, $(2\,\Gamma (1+k/2)/\pi ^{k/2})W_{k-1}$
coincides with $MW$ when $k=3$. \ This is as far as our analysis has gone.

\section{Trihedron}

Consider now three unit spheres in $\mathbb{R}^{3}$ that pass through each
other's centers.\ \ The region enclosed by all spheres is called a \textbf{%
spherical trihedron}, using language of \cite{G0, G3, G1, G2}. Starting with
spheres%
\[
\left( x-\dfrac{1}{\sqrt{3}}\right) ^{2}+y^{2}+z^{2}=1, 
\]%
\[
\begin{array}{ccc}
\left( x+\dfrac{1}{2\sqrt{3}}\right) ^{2}+\left( y-\dfrac{1}{2}\right)
^{2}+z^{2}=1, &  & \left( x+\dfrac{1}{2\sqrt{3}}\right) ^{2}+\left( y+\dfrac{%
1}{2}\right) ^{2}+z^{2}=1%
\end{array}%
\]%
and defining%
\[
\begin{array}{ccccc}
f(x,y)=\sqrt{1-\left( x-\dfrac{1}{\sqrt{3}}\right) ^{2}-y^{2}}, &  & a(y)=%
\dfrac{1}{\sqrt{3}}-\sqrt{1-y^{2}}, &  & c(y)=-\dfrac{1}{\sqrt{3}}y%
\end{array}%
\]%
we obtain%
\[
VL=12\dint\limits_{0}^{1/2}\;\;\dint\limits_{a(y)}^{c(y)}f(x,y)dx\,dy=\frac{1%
}{12}\left[ 2\sqrt{2}+24\pi -57\func{arcsec}(3)\right] , 
\]%
\[
AR=12\dint\limits_{0}^{1/2}\;\;\dint\limits_{a(y)}^{c(y)}\sqrt{%
1+f_{x}^{2}+f_{y}^{2}}\,dx\,dy=6\left[ \pi -2\func{arcsec}(3)\right] . 
\]%
Each of the three edges is a circular arc with radius $\sqrt{3}/2$. Focus on
the arc that lies entirely in the $xz$-plane. From%
\[
\left( x+\dfrac{1}{2\sqrt{3}}\right) ^{2}+z^{2}=1-\dfrac{1}{4}=\left( \dfrac{%
\sqrt{3}}{2}\right) ^{2} 
\]%
we deduce that the circle center is $(-1/(2\sqrt{3}),0,0)$ and hence the
subtended angle is%
\[
2\arccos \left( \frac{1/(2\sqrt{3})}{\sqrt{3}/2}\right) =2\func{arcsec}(3) 
\]%
by symmetry. As before, we have%
\begin{eqnarray*}
MW &=&\frac{1}{2\pi }AR+\frac{3}{4\pi }\frac{\pi }{3}\left( 2\func{arcsec}%
(3)\right) \dfrac{\sqrt{3}}{2}=\frac{12\pi -\left( 24-\sqrt{3}\pi \right) 
\func{arcsec}(3)}{4\pi } \\
&=&1.182061751038757...\text{.}
\end{eqnarray*}%
Incidently, the vertex-to-vertex distance $\lambda $ here is $2\sqrt{2/3}$
and $VL/\lambda ^{3}=0.154...$, consistent with Figure 1 in \cite{G2}. \ See 
\cite{He} for extended volume results.

\section{Tetrahedron}

Consider finally four unit spheres in $\mathbb{R}^{3}$ that pass through
each other's centers.\ \ The region enclosed by all spheres is called a 
\textbf{spherical tetrahedron} or \textbf{Reuleaux tetrahedron}. Our
analysis of this solid is complicated due to a narrow sliver in the
half-space $z<0$. Starting with spheres%
\[
\begin{array}{ccc}
\left( x-\dfrac{1}{\sqrt{3}}\right) ^{2}+y^{2}+z^{2}=1, &  & \left( x+\dfrac{%
1}{2\sqrt{3}}\right) ^{2}+\left( y-\dfrac{1}{2}\right) ^{2}+z^{2}=1,%
\end{array}%
\]%
\[
\begin{array}{ccc}
\left( x+\dfrac{1}{2\sqrt{3}}\right) ^{2}+\left( y+\dfrac{1}{2}\right)
^{2}+z^{2}=1, &  & x^{2}+y^{2}+\left( z-\sqrt{\dfrac{2}{3}}\right) ^{2}=1%
\end{array}%
\]%
and defining%
\[
\begin{array}{ccc}
f(x,y)=\sqrt{1-\left( x-\dfrac{1}{\sqrt{3}}\right) ^{2}-y^{2}}, &  & g(x,y)=%
\sqrt{\dfrac{2}{3}}-\sqrt{1-x^{2}-y^{2}},%
\end{array}%
\]%
\[
\begin{array}{ccccc}
a(y)=\dfrac{1}{\sqrt{3}}-\sqrt{1-y^{2}}, &  & b(y)=\dfrac{1}{\sqrt{6}}\left( 
\dfrac{1}{\sqrt{2}}-\sqrt{3-4y^{2}}\right) , &  & c(y)=-\dfrac{1}{\sqrt{3}}y%
\end{array}%
\]%
we obtain \cite{Hb, Ws}%
\begin{eqnarray*}
VL
&=&12\dint\limits_{0}^{1/2}\;\;\dint\limits_{a(y)}^{b(y)}f(x,y)dx\,dy+6\dint%
\limits_{0}^{1/2}\;\;\dint\limits_{b(y)}^{c(y)}\left[ f(x,y)-g(x,y)\right]
dx\,dy \\
&=&\frac{1}{12}\left[ 3\sqrt{2}+32\pi -81\func{arcsec}(3)\right] ,
\end{eqnarray*}%
\begin{eqnarray*}
AR &=&12\dint\limits_{0}^{1/2}\;\;\dint\limits_{a(y)}^{b(y)}\sqrt{%
1+f_{x}^{2}+f_{y}^{2}}\,dx\,dy+6\dint\limits_{0}^{1/2}\;\;\dint%
\limits_{b(y)}^{c(y)}\left( \sqrt{1+f_{x}^{2}+f_{y}^{2}}+\sqrt{%
1+g_{x}^{2}+g_{y}^{2}}\right) \,dx\,dy \\
&=&2\left[ 4\pi -9\func{arcsec}(3)\right] .
\end{eqnarray*}%
Each of the six edges is a circular arc with radius $\sqrt{3}/2$. Focus
again on the arc that lies entirely in the $xz$-plane. Previously it started
at $(0,0,\sqrt{2/3})$ and ended at $(0,0,-\sqrt{2/3})$; now it ends at $(1/%
\sqrt{3},0,0)$ where the arc meets the bottom face. Hence the subtended
angle is half its previous value. We consequently have%
\begin{eqnarray*}
MW &=&\frac{1}{2\pi }AR+\frac{6}{4\pi }\frac{\pi }{3}(\func{arcsec}(3))%
\dfrac{\sqrt{3}}{2}=\frac{16\pi -\left( 36-\sqrt{3}\pi \right) \func{arcsec}%
(3)}{4\pi } \\
&=&1.006582094946935...\text{.}
\end{eqnarray*}%
which falls between the minimum width $1$ and the maximum width $\sqrt{3}-1/%
\sqrt{2}=1.0249...$, as expected \cite{KW}.

Incidently, the vertex-to-vertex distance $\lambda $ here is $1$ and $%
VL/\lambda ^{3}=0.422...$, consistent with Figure 1 in \cite{G2}. \ See \cite%
{He} for extended volume results.

The \textbf{Meissner tetrahedron}, obtained by rounding off three edges of
the Reuleaux tetrahedron to force a constant width, possesses measures \cite%
{KW, By, Hr, CG}%
\[
VL^{\prime }=\frac{1}{12}\left[ 8-3\sqrt{3}\func{arcsec}(3)\right] \pi <VL, 
\]%
\[
AR^{\prime }=\frac{1}{2}\left[ 4-\sqrt{3}\func{arcsec}(3)\right] \pi <AR 
\]%
and, of course, $MW^{\prime }=1<MW$.

\section{Miscellanea}

The symmetric lens arises as a solution of certain geometric optimization
problems \cite{U1, U2, U3}. \ We report on three other such
\textquotedblleft important\textquotedblright\ solids here. \ In the same
papers, the \textbf{right circular cylinder with hemispherical ends} is
featured. \ Assuming the radius is $1$ and the cylinder length is $\ell $,
we easily have \cite{JP}%
\[
VL=\left( \ell +\frac{4}{3}\right) \pi , 
\]%
\[
AR=2\left( \ell +2\right) \pi , 
\]%
\[
MW=\frac{1}{2}\left( \ell +4\right) . 
\]%
The measures of a hemisphere, as an aside, are given in standard tables \cite%
{SKM, Sa} (although not spherical caps in general).

The \textbf{symmetric segment} or \textbf{spherical slice }appears in \cite%
{V1, V2, V3, V4}. This solid is obtained by removing two
diametrically-opposed spherical caps from the unit ball, each of angular
radius $\varphi $. \ Clearly \cite{Hw1} 
\[
VL=\frac{2\pi }{3}\left( 2+\sin (\varphi )^{2}\right) \cos (\varphi ), 
\]%
\[
AR=2\pi \left( 2\cos (\varphi )+\sin (\varphi )^{2}\right) , 
\]%
\[
MW=2\cos (\varphi )+\varphi \,\sin (\varphi ). 
\]%
Figure 5e of \cite{Hw1} is unambiguous (unlike Figure 5f, as discussed
earlier). \ \ We are, however, unable to find agreement with the
quermassintegrals in \cite{Hw2} (formula 68) when $k=3$.

The \textbf{cap body of a ball} appears in the same papers, as well as in 
\cite{W1, W2, W3, W4, W5, W6}. \ This solid is the convex hull of the unit
ball with a line segment passing symmetrically through its center. Another
name for this is\ \textbf{1-tangential body}$.$ With an interpretation of $%
\varphi $ identical to above, we find that \cite{Hw1}%
\[
VL=\frac{2\pi }{3}\frac{1+\cos (\varphi )^{2}}{\cos (\varphi )}, 
\]%
\[
AR=2\pi \frac{1+\cos (\varphi )^{2}}{\cos (\varphi )}, 
\]%
\[
MW=\frac{1+\cos (\varphi )^{2}}{\cos (\varphi )}. 
\]%
Figure 5g of \cite{Hw1} is unambiguous. \ Again, we are unable to find
agreement with the quermassintegrals in \cite{Hw2} (formula 70) when $k=3$.

A less important example (evidently) might be called the \textquotedblleft
ring body of a ball\textquotedblright . \ This solid is the convex hull of
the unit ball with a circle suspended symmetrically above its equator. \ Its
measures are given in \cite{Hw1} -- see Figure 5h -- but since it does not
appear elsewhere in the literature, we omit further discussion.

In closing, here is an unanswered question. The region enclosed by the six
spheres: 
\[
\begin{array}{ccc}
\left( x-\dfrac{1}{\sqrt{2}}\right) ^{2}+y^{2}+z^{2}=1, &  & \left( x+\dfrac{%
1}{\sqrt{2}}\right) ^{2}+y^{2}+z^{2}=1,%
\end{array}%
\]%
\[
\begin{array}{ccc}
x^{2}+\left( y-\dfrac{1}{\sqrt{2}}\right) ^{2}+z^{2}=1, &  & x^{2}+\left( y+%
\dfrac{1}{\sqrt{2}}\right) ^{2}+z^{2}=1,%
\end{array}%
\]%
\[
\begin{array}{ccc}
x^{2}+y^{2}+\left( z-\dfrac{1}{\sqrt{2}}\right) ^{2}=1, &  & 
x^{2}+y^{2}+\left( z+\dfrac{1}{\sqrt{2}}\right) ^{2}=1%
\end{array}%
\]%
is called a \textbf{spherical hexahedron (cube)}. What are exact expressions
for $VL$, $AR$ and $MW$? \ The dihedral angle at any edge is $\pi /3$. While
the spheres are not mutually equidistant, we can still define $\lambda $ to
be the \textit{adjacent} vertex-to-vertex distance, which here is $2/\sqrt{3}
$. \ Numerical work gives $VL/\lambda ^{3}=1.508...$, consistent with Figure
1 in \cite{G2}. \ A similar question can be asked about the spherical
dodecahedron, for which $VL/\lambda ^{3}$ was given in \cite{G2} to be
approximately $7.86$.

\section{Acknowledgements}

Wouter Meeussen's package ConvexHull3D.m was helpful to me in preparing this
paper \cite{Msn}. \ He kindly extended the software functionality at my
request. \ I\ would be grateful for assistance in expanding my bibliography:
surely I have missed more than a few documents on measures of intersections
of balls!

\begin{tabular}{lll}
& Steven R. Finch &  \\ 
& Dept. of Statistics &  \\ 
& Harvard University &  \\ 
& Cambridge, MA, USA &  \\ 
& \textit{steven\_finch@harvard.edu} & 
\end{tabular}


\begin{thebibliography}{99}
\bibitem{Fi1} S. R. Finch, Convex hull of two orthogonal disks,
http://arxiv.org/abs/1211.4514.

\bibitem{Fi2} S. R. Finch, Oblique circular cones and cylinders,
http://arxiv.org/abs/1212.5946.

\bibitem{Rt} G.-C. Rota, Mathematical snapshots, unpublished note (1997),
http://www-groups.dcs.st-and.ac.uk/\symbol{126}history/Extras/rota.pdf.

\bibitem{G0} M. E. Glicksman, Analysis of 3-D network structures, \textit{%
Philosophical Magazine} 85 (2005) 3--31.

\bibitem{G3} M. E. Glicksman, P. R. Rios and D. J. Lewis, Regular $N$-hedra:
A topological approach for analyzing three-dimensional textured
polycrystals, \textit{Acta Materialia} 55 (2007) 4167--4180.

\bibitem{G1} M. E. Glicksman, P. R. Rios and D. J. Lewis, Mean width and
caliper characteristics of network polyhedra, \textit{Philosophical Magazine}
89 (2009) 389--403.

\bibitem{G2} M. E. Glicksman, P. R. Rios and D. J. Lewis, Linear measures
for polyhedral networks, \textit{Internat. J. Materials Res.} 100 (2009)
536-542.

\bibitem{S1} H. L. Weissberg and S. Prager, Viscous flow through porus
media. II. Approximate three-point correlation function, \textit{Physics of
Fluids} 5 (1962) 1390--1392.

\bibitem{S2} J. S. Rowlinson, The triplet distribution function in a fluid
of hard spheres, \textit{Molecular Physics} 6 (1963) 517--524.

\bibitem{S3} M. J. D. Powell, The volume internal to three intersecting hard
spheres, \textit{Molecular Physics} 7 (1964) 591--592.

\bibitem{S4} F. H. Ree, R. N. Keeler, and S. L. McCarthy, Radial
distribution function of hard spheres, \textit{J. Chem. Physics} 44 (1966)
3407--3425.

\bibitem{S5} R. Lustig, Surface and volume of three, four, six and twelve
hard fused spheres, \textit{Molecular Physics} 55 (1985) 305--317.

\bibitem{S6} R. Lustig, Geometry of four hard fused spheres in an arbitrary
spatial configuration, \textit{Molecular Physics} 59 (1986) 195--207.

\bibitem{S7} K. D. Gibson and H. A. Scheraga, Volume of the intersection of
three spheres of unequal size: A simplified formula, \textit{J. Physical
Chem. }91 (1987) 4121-4122; addenda 91 (1987) 6326.

\bibitem{S8} K. D. Gibson and H. A. Scheraga, Surface area of the
intersection of three spheres with unequal radii: A simplified analytical
formula, \textit{Molecular Physics} 64 (1988) 641--644.

\bibitem{S9} L. R. Dodd and D. N. Theodorou, Analytical treatment of the
volume and surface area of molecules formed by an arbitrary collection of
unequal spheres intersected by planes, \textit{Molecular Physics} 72 (1991)
1313--1345.

\bibitem{S0} L. S. Chkhartishvili, Volume of the domain of intersection of
three spheres (in Russian), \textit{Mat. Zametki} 69 (2001) 466--476; Engl.
transl. in \textit{Math. Notes} 69 (2001) 421--428; MR1846843 (2002d:52004).

\bibitem{He} A. Helte, Fourth-order bounds on the effective conductivity for
a system of fully penetrable spheres, \textit{Proc. Royal Soc. London A} 445
(1994) 247--256.

\bibitem{Hw1} H. Hadwiger, \textit{Altes und Neues \"{u}ber konvexe K\"{o}%
rper}, Birkh\"{a}user Verlag, 1955, pp 35--37; MR0073220 (17,401e).

\bibitem{Hw2} H. Hadwiger, \textit{Vorlesungen \"{u}ber Inhalt, Oberfl\"{a}%
che und Isoperimetrie}, Springer-Verlag, 1957, pp. 215--221; MR0102775 (21
\#1561).

\bibitem{WL} S. R. Wilson, C. M. Hefferan, S. F. Li, J. Lind, R. M. Suter
and A. D. Rollett, Microstructural characterization and evolution in 3D, 
\textit{Challenges in Materials Science and Possibilities in 3D and 4D
Characterization Techniques}, Proc. 31$^{\text{st}}$ Ris\o\ Internat. Symp.
on Materials Science, ed. N. Hansen, D. Juul Jensen, S. F. Nielsen, H. F.
Poulsen, and B. Ralph, Ris\o\ National Laboratory for Sustainable Energy,
Technical University of Denmark, 2010, pp. 201--217;
http://mimp.mems.cmu.edu/publications/.

\bibitem{Oe} M. Oettel, H. Hansen-Goos, P. Bryk and R. Roth, Depletion
interaction of two spheres - Full density functional theory vs. morphometric
results, \textit{Europhysics Letters} 85 (2009) 36003.

\bibitem{T1} Wikipedia contributors, Spherical cap, \textit{Wikipedia, The
Free Encyclopedia}, 22 July 2009,
http://en.wikipedia.org/wiki/Spherical\_cap.

\bibitem{T2} S. Li, Concise formulas for the area and volume of a
hyperspherical cap, \textit{Asian J. Math. Stat.} 4 (2011) 66--70; MR2813331.

\bibitem{T3} J. S. Brauchart, D. P. Hardin and E. B. Saff, The next-order
term for optimal Riesz and logarithmic energy asymptotics on the sphere, 
\textit{Recent Advances in Orthogonal Polynomials, Special Functions, and
their Applications}, Proc. 11$^{\text{th}}$ Internat. Symp. (OPSFA'11), Legan%
\'{e}s, ed. J. Arves\'{u} and G. L\'{o}pez Lagomasino, Amer. Math. Soc.,
2012, pp. 31--61; MR2964138; http://arxiv.org/abs/1202.4037.

\bibitem{T4} G. D. Anderson, M. K. Vamanamurthy and M. K. Vuorinen, \textit{%
Conformal Invariants, Inequalities, and Quasiconformal Maps}, John Wiley \&
Sons, 1997, pp. 38--41, 163--164; MR1462077 (98h:30033).

\bibitem{Hb} B. Harbourne, Volume and surface area of the spherical
tetrahedron (AKA Reuleaux tetrahedron) by geometrical methods,
http://www.math.unl.edu/\symbol{126}bharbourne1/ST/sphericaltetrahedron.html.

\bibitem{Ws} E. W. Weisstein, Reuleaux tetrahedron, \textit{MathWorld - A
Wolfram Web Resource}, http://mathworld.wolfram.com/ReuleauxTetrahedron.html.

\bibitem{KW} B. Kawohl and C. Weber, Meissner's mysterious bodies,\textit{\
Math. Intelligencer} 33 (2011) 94--101; MR2844102 (2012j:52006).

\bibitem{By} T. Bayen, T. Lachand-Robert and \'{E}. Oudet, Analytic
parametrization of three-dimensional bodies of constant width, \textit{Arch.
Ration. Mech. Anal.} 186 (2007) 225--249; MR2342202 (2008f:52005).

\bibitem{Hr} E. M. Harrell, Calculations for convex bodies. Example: the
rotated Reuleaux triangle, unpublished note (2001).

\bibitem{CG} G. D. Chakerian and H. Groemer, Convex bodies of constant
width, \textit{Convexity and its Applications}, ed. P. M. Gruber and J. M.
Wills, Birkh\"{a}user, Basel, 1983, pp. 49--96; MR0731106 (85f:52001).

\bibitem{U1} J. Favard, Sur quelques probl\`{e}mes de couvercles, \textit{%
Colloque de G\'{e}om\'{e}trie Diff\'{e}rentielle, Louvain}, Georges Thone,
Masson \& Cie, 1951, pp. 37--49; MR0050296 (14,309d). \ 

\bibitem{U2} H. Hadwiger, Elementare Ermittlung extremaler Rotationsk\"{o}%
rper, \textit{Revista Mat. Hisp.-Amer.} 9 (1949) 59--70; MR0035040 (11,680g).

\bibitem{U3} H. Hadwiger, Neue Ungleichungen f\"{u}r konvexe Rotationsk\"{o}%
rper, Math. Annalen 122 (1950) 175--180; MR0039290 (12,526c).

\bibitem{JP} K. M. Jansons and C. G. Phillips, On the application of
geometric probability theory to polymer networks and suspensions. I, \textit{%
J. Colloid and Interface Science} 137 (1990) 75--93.

\bibitem{SKM} D. Stoyan, W. S. Kendall and J. Mecke, \textit{Stochastic
Geometry and its Applications}, Wiley, 1987, pp. 11--19; MR0895588
(88j:60034a).

\bibitem{Sa} L. A. Santal\'{o}, \textit{Integral Geometry and Geometric
Probability}, Addison-Wesley, 1976, pp. 226--230; MR0433364 (55 \#6340).

\bibitem{V1} H. Hadwiger, P. Glur and H. Bieri, Die symmetrische Kugelzone
als extremaler Rotationsk\"{o}rper, \textit{Experienta} 4 (1948) 304--305;
MR0026348 (10,141d).

\bibitem{V2} H. Bieiri, Mitteilung zum Problem eines konvexen Extremalk\"{o}%
rpers, \textit{Arch. Math. (Basel)} 1 (1949) 462--463; MR0031279 (11,127c).

\bibitem{V3} H. Hadwiger, Beweis einer Extremaleigenschaft der symmetrischen
Kugelzone, \textit{Portugaliae Math.} 7 (1948) 73--85; MR0028597 (10,471a).

\bibitem{V4} M. A. Hern\'{a}ndez Cifre, G. Salinas and S. Segura Gomis, Two
optimization problems for convex bodies in the $n$-dimensional space, 
\textit{Beitr\"{a}ge Algebra Geom.} 45 (2004) 549--555; MR2093025
(2005e:52008).

\bibitem{W1} G. Bol, Beweis einer Vermutung von H. Minkowski, \textit{%
Abhandlungen aus dem Mathematischen Seminar der Hansischen Universit\"{a}t}
15 (1943) 37--56; MR0015824 (7,474f).

\bibitem{W2} J. R. Sangwine-Yager, The missing boundary of the Blaschke
diagram, \textit{Amer. Math. Monthly} 96 (1989) 233--237; MR0991869
(90a:52024).

\bibitem{W3} J. R. Sangwine-Yager, Stability for a cap body inequality, 
\textit{Geom. Dedicata} 38 (1991) 347--355; MR1112672 (92j:52012).

\bibitem{W4} V. A. Zalgaller, A family of extremal spindle-shaped bodies (in
Russian), \textit{Algebra i Analiz} 5 (1993) 200--214; Engl. transl. in 
\textit{St. Petersburg Math. J.} 5 (1994) 177--188; MR1220496 (94f:52009).

\bibitem{W5} S. Campi, Three-dimensional Bonnesen type inequalities, \textit{%
Le Matematiche (Catania)} 60 (2005) 425--431; MR2257044 (2007g:52008).

\bibitem{W6} S. Campi and P. Gronchi, A Favard-type problem for 3d convex
bodies, \textit{Bull. Lond. Math. Soc.} 40 (2008) 604--612; MR2441132
(2009h:52021).

\bibitem{Msn} W. Meeussen, Various Mathematica files,
http://users.telenet.be/Wouter.Meeussen/.

\bibitem{Fi3} S. R. Finch, Various preprints about mean width and intrinsic
volumes, http://www.people.fas.harvard.edu/\symbol{126}sfinch/csolve/.
\end{thebibliography}
\end{document}